\documentclass[12pt]{article}
\usepackage{cite}
\usepackage{amssymb}
\usepackage{amsmath}
\usepackage{graphics}
\numberwithin{equation}{section}

\newcommand{\eh}{\hfill}\newlength{\sperr}

\newenvironment{proof}{{\settowidth{\sperr}{\bf\rm
Proof}%
\par\addvspace{0.3cm}\noindent\parbox[t]{1.3\sperr}
{\textit{ P\eh r\eh o\eh o\eh f\eh .}}%
}}{\nopagebreak\mbox{}
$\Box$\par\addvspace{0.3cm}}

\def\ve{\varepsilon}

\def\BC{{\mathbb C}}

\newtheorem{Pa}{Paper}[section]
\newtheorem{Tm}[Pa]{{\bf Theorem}}
\newtheorem{La}[Pa]{{\bf Lemma}}
\newtheorem{Cy}[Pa]{{\bf Corollary}}
\newtheorem{Rk}[Pa]{{\bf Remark}}

\newtheorem{Ee}[Pa]{{\bf Example}}
\newtheorem{Dn}[Pa]{{\bf Definition}}
\newtheorem{Pn}[Pa]{{\bf Proposition}}
\numberwithin{equation}{section}

\title{$(S+N)$-triangular operators: spectral properties and important examples}

\author{Lev Sakhnovich}

\date{}

\begin{document}

\maketitle

\vspace{0.5em}

\emph{99 Cove ave., Milford, CT, 06461, USA; \\
 E-mail: lsakhnovich@gmail.com}

 \vspace{0.5em}

\textbf{Abstract.} We introduce a notion of $(S+N)$-triangular  operators in the Hilbert space
using some basic ideas from  triangular representation theory.
Our notion generalizes the well-known notion of the spectral operators so that many properties
of the  $(S+N)$-triangular  operators coincide with those of spectral operators. At the same time
we show that wide classes of operators are $(S+N)$-triangular.

\textbf{Mathematics Subject Classification 2010:} Primary 47A45; \\
Secondary 47A15, 47A46.

\textbf{Keywords.} Maximal chain of orthogonal projections, invariant subspace, 
triangular
representation, spectral operator.

\section{Main notions}1. Let us consider bounded linear operators
acting in the separable Hilbert  space $H$. The next three definitions introduce the notions
of a {\it quasinilpotent operator} and of a {\it scalar type operator}. Recall that the notion of the scalar type operator
is introduced  in terms of the {\it resolution of identity} which, in turn,
is introduced  via the notion of a spectral measure.
\begin{Dn}\label{Definition 1.1} A bounded operator $N$ is a quasinilpotent operator if and only if its spectrum
$\sigma(N)$ is concentrated at zero.\end{Dn}
\begin{Dn}\label{Definition 1.2} A bounded operator $S$ is said to be of scalar type if
\begin{equation}S=\int_{\BC}\lambda{E(d\lambda)},
\label{1.1}\end{equation}
where $\BC$ is the complex plane and $E$ is the resolution of the identity for $S$.\end{Dn}
\begin{Dn}\label{Definition 1.3} Spectral measure $E$ is the resolution of the identity for the operator $S$
if
\begin{equation}E(\delta)S=SE(\delta),\quad \sigma(S_{\delta}){\subseteq}\delta.\label{1.2}\end{equation}
\end{Dn}
Here $\delta$ stands for Borel sets and the symbol $S_{\delta}$ denotes the restriction of $S$ to the subspace
$H(\delta)=E(\delta)H$.

Using the notions of  quasinilpotent and  scalar type operators, N. Dunford  \cite{Dun1, Dun2} introduced the notion of a spectral operator and investigated the class of such operators. 
This theory was actively developed in numerous papers and books (see, e.g., \cite{A, AA, AR, BaMe, CoFo, DunSc, GeT, MMN} and references therein).
\begin{Dn}\label{Definition 1.4} A bounded operator $T$ is a spectral operator if it is a sum
\begin{equation}T=S+N \label{1.3}\end{equation}
of a bounded scalar type operator $S$ and of a quasinilpotent operator $N$ which commutes with $S$. \end{Dn}
It is proved (see \cite{DunSc}) that the decomposition \eqref{1.3} is unique.

2. However, the requirement that $S$ and $N$ commute is very strong. Various important operators  do not have the "spectral" property.
In the present paper we use some ideas of the spectral operators theory  \cite{Dun1}, \cite{Dun2} and triangular representation theory \cite{Sakh1, Sakh2, Sakh3}.
We again consider the representation of the operator $T$ in the form \eqref{1.3},
 where $S$ is a bounded scalar type operator and $N$ is a quasinilpotent operator.
However,  we do not assume that the operator $N$ commutes  with the operator $S$.
Let us formulate our assumption.

We assume that there exists such a maximal chain of  orthogonal projections $E_t$ $(a \leq t \leq b)$
that
\begin{equation}\label{1.4} E_t{S}={S}E_t ,\quad E_t{N}={N}E_t ,\quad a<t<b.
\end{equation}
Let us remind the notion of the maximal chain.
\begin{Dn}\label{Definition 1.5}The chain of
the orthogonal projections $E_t$ $(a \leq t \leq b)$ is maximal if $\dim(E_{a+0}-E_{a}){\leq}1$,  $\, \dim(E_{b}-E_{b-0}){\leq}1$ and
\begin{equation}E_{a}H=0,\quad E_{b}H=H, \quad \dim(E_{t+0}-E_{t-0}){\leq}1 \quad
(a<t<b).\label{1.5}\end{equation}
\end{Dn}
\begin{Rk}\label{Remarkd1} Although $E_t$ differs from a spectral measure, some maximal chain $E_t$ corresponding to a scalar type
operator $S$ always exists, that is, the first equality in \eqref{1.4} holds.
Moreover, some maximal chain $E_t$ always exists for a compact operator $N$.  In \eqref{1.4} we assume that $S$ and $N$
have the same maximal chain $E_t$.
\end{Rk}
It follows from  \eqref{1.4} that
\begin{equation} E_t{T}E_t={T}E_t,\quad a<t<b.\label{1.6}\end{equation}
According to \eqref{1.6} the subspaces $H_t=E_t{H}$ are invariant subspaces of the operator $T$. Hence, formula \eqref{1.6} gives the triangular structure  of the operator $T$ (see \cite{Sakh1, Sakh2, Sakh3}).

Consider $(S-{\lambda}I)^{-1}$ and $(T-{\lambda}I)^{-1}$, where $I$ is the identity operator.\\
\textbf{Conditions on the quasinilpotent operator $N$}\label{Addition condition}:

$\bf{( A)}$  \emph{Let $\lambda$ be a regular point of the resolvent
$(S-{\lambda}I)^{-1}$.  Then the operator
$(S-{\lambda}I)^{-1}N$  is quasinilpotent.}

$\bf{(B)}$
\emph{Let $\lambda$ be a regular point of the resolvent
$(T-{\lambda}I)^{-1}$.  Then the operator
$(T-{\lambda}I)^{-1}N$  is quasinilpotent.}

\begin{Rk}\label{Remark 1.6}1. If $SN=NS$ $($i.e., if the operator $T$ is spectral$)$, then the conditions $\bf{(A)}$ and $\bf{(B)}$ on the operator $N$  are fulfilled.

2. If the operator $N$ is compact,  then the conditions $\bf{(A)}$ and $\bf{(B)}$  on the operator $N$ are fulfilled $($see \cite[Ch.2, Section 17]{Br}$)$.\end{Rk}

\begin{Dn}\label{Definition 1.7} We say that  the operator $T$ is  $(S+N)$-triangular $($where $S$ is a bounded scalar type operator and $N$ is a quasinilpotent operator$)$,
if relations \eqref{1.3} and \eqref{1.4} as well as conditions $\bf{(A)}$ and $\bf{(B)}$  are fulfilled.
\end{Dn}
It is easy to see that all the spectral operators are (S+N)-triangular as well.
\begin{Rk}\label{Remark 1.8}   All  the operators  in the finite  dimensional  spaces   are $(S~+~N)$-triangular.\end{Rk}
In many important cases the corresponding operator $N$ is compact.
\begin{Dn}\label{Definition 12} We say that the operator T is $(S+V)$-triangular if  $T$ is
$(S+N)$-triangular  and the corresponding operator $N$ is compact.\end{Dn}
It is evident that the set of (S+V)-triangular   operators is a subclass of the class of
$(S+N)$-triangular operators. The general theory of $(S+N)$-operators and $(S+V)$-operators is developed in section 2.
In sections 3 and 4, we apply the general theory  of the $(S+V)$-triangular operators
to the operators of the form
\begin{equation}Tf=\phi(x)f(x)+\int_{0}^{x}v(x,t)f(t)dt.\label{1.15}\end{equation}
Spectral mapping theorems,  conditions under which $T$ is decomposable (see Definition \ref{Definition 3.4} for the notion
of a decomposable operator) and
estimates for the resolvent $R_{\lambda}(T)=({\lambda}I-T)^{-1}$ are given.
In the last section 5 we study triangular integro-differential operators with difference kernels
\begin{equation}Tf=\frac{d}{dx}\int_{0}^{x}s(x-t)f(t)dt,\quad f{\in}H=L^{2}(0,\omega).\label{1.16}\end{equation}
In particular, we show that in the general case situation the operators of the form \eqref{1.16} are not (S+V)-triangular.

\section{(S+N)-triangular operators: spectral mapping theorems}
We need  the following well-known assertion (see for example\cite{DunSch}, Theorem VII.3.11).
\begin{Tm}\label{Theorem 2.1}Let the operator $T$ be $(S+N)$-triangular. Then
\begin{equation}\sigma(T)=\sigma(S).\label{2.1}\end{equation}
\end{Tm}
\begin{proof}
 If $\lambda{\notin}\sigma(S)$, then we have
\begin{equation}(T-\lambda{I})=(S-\lambda{I})\big(I+(S-\lambda{I})^{-1}N\big).\label{2.2}\end{equation}
According to the condition  $\bf{(A)}$    on the operator $\,N$, the operator $(S~-~\lambda{I})^{-1}N$ is quasinilpotent.
In view of \eqref{2.2} the operator $(T-\lambda{I})$ is invertible and $\lambda{\notin}\sigma(T)$. In the same way we prove that  
the relation $\lambda{\notin}\sigma(S)$ follows from the relation
$\lambda{\notin}\sigma(T)$. The theorem is proved.
\end{proof}
\begin{Cy}\label{Corollary 2.2} Let the operator $T$ be $(S+V)$-triangular.
Then, for any polynomial $p(z)$, the  operator $p(T)$ is an $(S+V)$-triangular operator  and
\begin{equation}p(T)=p(S)+V_{p}\, ,\quad \sigma(p(T))=\sigma(p(S))=p(\sigma(S))=p(\sigma(T)),\label{2.3}\end{equation}
where the operator $V_{p}$ is quasinilpotent and compact.\end{Cy}
The following spectral mapping theorem is valid.
\begin{Tm}\label{Theorem 2.3}Let the operator $T$ be $(S+V)$-triangular.  Then, for any function $f(z)$ which
is analytic in some neighborhood  of $\sigma(T)$, we have
\begin{equation} \sigma(f(T))=f(\sigma(T)).\label{2.4}\end{equation}
\end{Tm}
\begin{proof} Using  relation \eqref{2.3} we see that relations \eqref{2.4} hold when
$f(z)$ is a polynomial (i.e., $f(z)=p(z)$). Let us consider the case
$$f(z)=r(z)=~p(z)/q(z),$$
where $p(z)$ and $q(z)$ are polynomials
(i.e., $f(z)$ is a rational function) and $q(z)$ has no zeros
in $\sigma(T)$. In this case, we have
\begin{equation}r(T)-{\lambda}I=\big(p(T)-{\lambda}q(T)\big)\big(q(T)\big)^{-1}.\label{2.5}\end{equation}
Taking into account Corollary \ref{Corollary 2.2} and equality \eqref{2.5} we derive that the following assertions are equivalent:\\
(i) $\lambda{\notin}\sigma(r(T))$;\\
(ii) the polynomial $p(z)-{\lambda}q(z)$ has no zeros in $\sigma(T)$;\\
(iii) the rational function $\big(p(z)-{\lambda}q(z)\big)\big(q(z)\big)^{-1}=r(z)-\lambda$ has no zeros in $\sigma(T)$.
From the assertions (i) and (iii) we have
\begin{equation} \sigma(r(T))=r(\sigma(T)).\label{2.6}\end{equation}
According to Runge's theorem (see, e.g.,  \cite{Mark}), if $f(z)$ is  analytic in a neighborhood of  $\sigma(T)$,  then  there is
a sequence $r_{n}(z)$  of rational functions with no poles in $\sigma(T)$  such that
$r_{n}(z)$ tends to $f(z)$ uniformly on $\sigma(T)$.
Thus, we obtain
\begin{equation}||r_{n}(T)-f(T)||{\to}0,\quad n{\to}\infty.\label{2.7}\end{equation}
Here we used the well-known formula
\begin{equation}f(T)=\frac{1}{2{\pi}i}\int_{C}f(\lambda)({\lambda}I-T)^{-1}d\lambda,
\label{2.8}
\end{equation}
where $C$  is an admissible contour surrounding the spectrum of T and
$f(z)$ is an analytic function in a neighborhood of $\sigma(T)$ . Relations \eqref{2.6} and
\eqref{2.7} imply the assertion of the theorem. \end{proof}
\begin{Pn}\label{Proposition 2.4}Let the operator $T$ be $(S+V)$-triangular and assume that $f(z)$
 is an analytic function in a neighborhood of $\sigma(T)$.
Then the  operator $f(T)$ is an $(S+V)$-triangular operator  and
\begin{equation}f(T)=f(S)+V_{f},  \label{2.9}\end{equation}
where the operator $V_{f}$ is quasinilpotent and compact.\end{Pn}
\begin{proof} It follows from \eqref{2.7} and relation
\begin{equation}||r_{n}(S)-f(S)||{\to}0,\quad n{\to}\infty\label{2.10}\end{equation}
that
\begin{equation}||V_{n}-V_{f}||{\to}0,\quad n{\to}\infty,\label{2.11}\end{equation}
where $V_{f}=f(T)-f(S)$.The operator $V_{n}$ is quasinilpotent and compact. Then in view of \eqref{2.11} the operator $V_{f}$ is quasinilpotent and compact too.
The proposition is proved. \end{proof}
According to \eqref{2.3},\eqref{2.4} and \eqref{2.9} we have:
\begin{Pn}\label{Proposition 2.5} Let the operator $T$ be $(S+V)$-triangular and assume that $f(z)$
 is an analytic function in a neighborhood of $\sigma(T)$. Then
 \begin{equation} \sigma(f(T))=\sigma(f(S))=f(\sigma(S))=f(\sigma(T)).\label{2.12}\end{equation}
 \end{Pn}
\section{$(S+N)$-triangular operators with  self-adjoint operators $S$}
1. Let us  consider  bounded operators $T$:
\begin{equation}T=Q+iB,\label{3.1}\end{equation}  acting in the Hilbert separable space $H$.
We assume that  the operator $Q$ is self-adjoint and the operator $B$ is self-adjoint and compact.
We need  the following definition  (see \cite{Mac, Gokr}).
\begin{Dn}\label{Definition 3.1} Denote by $\mathfrak{H}_{\omega}$ the Banach space of compact operators $B$ such that
\begin{equation}||B||_{\omega}:=\sum_{1}^{\infty}s_{n}(B)/(2n-1)<\infty,\label{3.2}\end{equation}
where the s-numbers $s_n(B)$ form the sequence of the non-decreasing eigenvalues of the operator $(B^{*}B)^{1/2}$.
\end{Dn}
Next, some conditions, under which  $T$ is $(S+V)$-triangular, are given.
\begin{Tm}\label{Theorem 3.2}Let the following conditions be fulfilled$:$\\
1.The operator $T$ has the form \eqref{3.1}.\\
2.The spectrum of the operator $T$ is real.\\
3.The corresponding operator $B$ belongs to the Banach space $\mathfrak{H}_{\omega}$.\\
Then the operator $T$ is $(S+V)$-triangular,
where the operator $S$ is self-adjoint and the operator $V$ is quasinilpotent and
compact.\end{Tm}
V.I. Macaev \cite{Mac} showed that the operator $T$ satisfying the conditions of the theorem has  non-trivial invariant subspaces. In 1957,
 we proved \cite{Sakh1} that in this case  there exists a maximal chain of orthogonal projections $E_t$ such that relation \eqref{1.6} is valid. 
 Later this result was repeated by M.S. Brodskii \cite{Br} and by J.R. Ringrose \cite{Ring} (in the years 1961 and 1962, respectively).

{\it Proof of Theorem \ref{Theorem 3.2}.}
We shall use \eqref{1.6} and  the following  statement \cite[Theorem 17.1]{Br}:\\
There exists a quasinilpotent and compact operator $V$ such that 
\begin{equation}
B=\frac{1}{2i}(V-V^{*}),\label{3.3}
\end{equation}
and the second equality in \eqref{1.4} (where $N=V$)   is valid.

It follows from \eqref{3.1} and \eqref{3.3} that the operator $S$ given by
\begin{equation} S=T-V\label{3.4}\end{equation}
is self-adjoint.  Moreover, in view of the second equality in \eqref{1.4}, of \eqref{1.6} and of \eqref{3.4} we have
\begin{equation}E_t{S}E_t=SE_t.\label{3.5}\end{equation}
Since the operator $S$ is self-adjoint, equality \eqref{3.5} implies the first equality in \eqref{1.4}.
Now, the assertion of the theorem follows directly from equality \eqref{3.4}.

2. The class of compact operators $B$ such that
\begin{equation}||B||_p:=\left(\sum_{n}s_{n}^{p}(B)\right)^{1/p}<\infty \quad (p{\geq}1)\label{3.6}\end{equation}
is denoted by $\mathfrak{H}_p$. Clearly, for $p\geq 1$ we have $\mathfrak{H}_p{\subset}\mathfrak{H}_{\omega}$. We note that subclasses of $(S+V)$-triangular operators from $\mathfrak{H}_p,$ $p\geq 2$  (without introducing any term
for these subclasses) already appeared in the literature.
More precisely,
Theorem 3.2 was proved in our paper \cite{Sakh2}  for the subcase $B\in \mathfrak{H}_2$ ($B$ is from the Hilbert-Schmidt class). Later  Theorem 3.2 was proved for the subcase 
$B \in \mathfrak{H}_p$ $(p>2)$, see \cite{Sc}.

3.  Let the operator $T$ be (S+V)-triangular and let the operator $S$ be self-adjoint.
 Using formula \eqref{2.2} we obtain
\begin{equation}R_{\lambda}(T)=(T-\lambda{I})^{-1}=\left(I+\sum_{n=1}^{\infty}
\frac{c_{n}(\lambda)}{(\Im{\lambda})^n}\right)(S-{\lambda}I)^{-1},\quad \Im{\lambda}{\ne}0,\label{3.7}\end{equation}
where
\begin{equation}c_{n}(\lambda)=\big(-(\Im{\lambda})V_{1}(\lambda)\big)^{n},\quad
V_{1}(\lambda)=(S-{\lambda}I)^{-1}V.\label{3.8}\end{equation}
Since $S$ is self-adjoint and  $V_{1}(\lambda)$ is quasinilpotent, we obtain the following proposition.
\begin{Pn}\label{Proposition 3.3} Let the operator $T$ be $(S+V)$-triangular, where $S$ is self-adjoint. Then
\begin{equation}||c_{n}(\lambda)||{\leq}M^{n},\quad
\lim_{n{\to}\infty} ||c_{n}(\lambda)||^{1/n}=0,\quad \Im{\lambda}{\ne}0.\label{3.9}\end{equation}\end{Pn}
Put $\lambda=x+iy$ and introduce the sequence of functions
\begin{equation}r_n(y)=\sup_{-\infty<x<\infty}||c_n(x+iy)||^{1/n}.
 \label{3.10}\end{equation}
Denote by  $\mathcal{N}(y)$ a number of values $r_n(y)$ which are greater than $|y|/2,\, |y|{\ne}0.$
It follows from \eqref{3.7}, \eqref{3.9} and \eqref{3.10}  that
\begin{equation}||R_{\lambda}(T)||{\leq}(C/|y|)(M/|y|)^{\mathcal{N}(y)}.
\label{3.11}
\end{equation}
Further we need some definitions from
\cite{CoFo}.
\begin{Dn}\label{Definition 3.5} A closed linear subspace $H_i$ is called a maximal spectral  subspace of $T$
if$:$ \\
1) The subspace $H_i$ is an invariant subspace of the operator $T$.\\
2) For any other   closed invariant subspace $\tilde{H}_i$ of   $T$,
such that \\ $\sigma(T\mid_{\tilde{H}_i}){\subseteq}\sigma(T\mid_{H_i})$, the relation $\tilde{H}_i{\subseteq}H_i$ is valid.
\end{Dn}
\begin{Dn}\label{Definition 3.4}An operator $T$ is called decomposable if for every finite open covering $\{G_i\}$ of $\sigma(T)$ there exists a system $\{H_i \}$ of maximal spectral  subspaces of $T$
such that $\sigma(T\mid_{H_i}){\subset}G_i$ for every $i$ and $H=\sum{H_i}$.\end{Dn}

\begin{Rk}\label{Remark 3.6}The class of decomposable operators from \cite{CoFo} is narrower than  the class
of S-operators from \cite{LyuMac}.
In the definition of an S-operator, the condition $H =\sum{H}_i$ is
replaced by a slightly weaker condition $H = \overline{\sum{H}_i}$. There are
examples satisfying this weaker condition which are not
decomposable in the sense of \cite{CoFo}, see [A].\end{Rk}
\begin{Dn}\label{Definition 3.7}Operator $T$ is called strongly decomposable if its restriction to an arbitrary maximal spectral subspace is decomposable too (see \cite{Ap}).\end{Dn}
\begin{Pn}\label{Proposition 3.8}If the inequality
\begin{equation}\int_{0}^{\ve} \ln({\mathcal{N}(y)})dy<\infty \label{3.12}\end{equation}
holds for some $\ve >0$, then the  operator $T$ is strongly decomposable.\end{Pn}
\begin{proof} According to \eqref{3.11} and \eqref{3.12} we have
\begin{equation}\int_{0}^{\ve} \ln \ln(\mathcal{M}(y))dy<\infty, \label{3.13}\end{equation}
where
\begin{equation}\mathcal{M}(y)=\sup_{\Im{\lambda}=y}||R_{\lambda}(T)||.\label{3.14}
\end{equation}
Inequality \eqref{3.13} is called Levinson's condition \cite{Lev, LyuMac}. According to \cite{Rad}, the assertion of the proposition follows from
Levinson's condition.
\end{proof}
\begin{Cy}\label{Corollary 3.10} Let the operator $T$ be $(S+V)$-triangular and let the operator $S$ be self-adjoint. If the inequalities
\begin{equation}\|c_{n}(\lambda)\|{\leq}\left(\frac{M}{ (\ln(n+1))^{\alpha}}\right)^{n}, \quad
\alpha>0\label{3.16} \end{equation} are valid, then
\begin{equation}\ln\mathcal{N}(y){<}({2M}/{|y|})^{1/\alpha},
\quad 0<|y|{\leq}\ve.\label{3.17}
\end{equation}\end{Cy}
Using Proposition \ref{Proposition 3.8} and Corollary \ref{Corollary 3.10}
we obtain the statement.
\begin{Cy}\label{Corollary 3.11}Let conditions of Corollary \ref{Corollary 3.10} be fulfilled
and let $\alpha>1$. Then the corresponding operator $T$ is strongly decomposable.\end{Cy}
\section{Examples}
\begin{Ee}\label{Example 4.1} Let the operator $T$ be $(S+V)$-triangular, where
\begin{equation}Vf=J^{\beta}f=
\frac{1}{\Gamma(\beta)}\int_{0}^{x}(x-t)^{\beta-1}f(t)dt,\quad \beta>0,\quad f{\in}H=L^{2}(0,\omega),\label{4.1}\end{equation}
and
\begin{equation} Sf=\phi(x)f(x).\label{4.2}\end{equation} 
$($Here the real valued function $\phi(x)$ is bounded and  $\Gamma(z)$ is the Euler gamma-function.$)$
\end{Ee}
In this case we have
\begin{equation}V^{n}f=J^{n\beta}f=
\frac{1}{\Gamma(n\beta)}\int_{0}^{x}(x-t)^{n\beta-1}f(t)dt,\quad \beta>0.\label{4.3}\end{equation}
It follows from \eqref{4.3} that
\begin{equation}||V^{n}||{\leq}\frac{1}{\Gamma(n\beta)}\int_{0}^{\omega}x^{n\beta-1}dx=
\frac{1}{\Gamma(n\beta+1)}{\omega}^{n\beta}.\label{4.4}\end{equation}
Using \eqref{3.8} and \eqref{4.4} we derive the inequality
\begin{equation}||c_{n}(\lambda)||{\leq}||V^{n}||{\leq}
\frac{1}{\Gamma(n\beta+1)}{\omega}^{n\beta}.\label{4.5}\end{equation}
Taking into account Stirling's formula, we obtain
\begin{equation}\Gamma(n\beta+1){\sim}\sqrt{2{\pi}n\beta}\left(\frac{n\beta }{e}\right)^{n\beta}.
\label{4.6}\end{equation}According to \eqref{4.5} and  \eqref{4.6},
 the function $\mathcal{N}(y)$   satisfies in this case  the relation
\begin{equation}\mathcal{N}(y){\leq}C|y|^{-1/\beta}. \label{4.7}\end{equation}
Then, in view of \eqref{3.11}, the following statement is valid.
\begin{Pn} \label{Proposition 4.2}Let the operators V and S be  defined by the relations \eqref{4.1} and \eqref{4.2}, respectively.
Then, the resolvent of the operator $T=S+V$  satisfies the inequalities
\begin{equation} \ln||R_{\lambda}(T)||{\leq}C|y|^{-1/\beta}|( \ln |y|)|,
\quad 0<|y|{\leq}\delta\label{4.8}
\end{equation}
for some $\delta >0$.
\end{Pn}
\begin{Cy}\label{Corollary 4.3}Let the operator $T$ have the form
\begin{equation}Tf=\phi(x)f(x)+\int_{0}^{x}v(x,t)f(t)dt,\quad f{\in}H=L^{2}(0,\omega),\label{4.9}\end{equation}where the real valued function $\phi(x)$
is bounded and  the kernel $v(x,t)$ is such that
\begin{equation}|v(x,t)|{\leq}M(x-t)^{\beta-1},\quad \beta>0.\label{4.10}\end{equation}
Then \\
1. The inequalities \eqref{4.7} and \eqref{4.8} are valid.\\
2.  The corresponding operator $T$ is strongly decomposable.\end{Cy}
\begin{Ee}\label{Example 4.4}Let the operator $T$ be $(S+V)$-triangular, where
\begin{equation}V_{\beta}f=
\frac{1}{\Gamma(\beta)}\int_{0}^{x}E_{\beta}(x-t)f(t)dt,\quad \beta>0,\quad f{\in}H=L^{2}(0,\omega),\label{4.11}\end{equation}
and
\begin{equation} Sf=\phi(x)f(x).\label{4.12}\end{equation} 
Here the real valued function $\phi(x)$ is bounded  and
\begin{equation}E_{\beta}(x)=\int_{0}^{\infty}\frac{1}{\Gamma(s)}e^{-Cs}s^{\beta-1}x^{s-1}ds,
\quad C=\overline{C}.\label{4.13}\end{equation}
\end{Ee}
In order to estimate $\mathcal{N}(y)$ we need some auxiliary results.
\begin{La}\label{Lemma 4.5} The function $E_{\beta}(x)$ is continuous in the domain $(0,\omega]$ and
\begin{equation}E_{\beta}(x)=\frac{\Gamma(\beta+1)}{x | \ln (x)|^{\beta+1}}\left[1+O\left(\frac{1}{ \ln (x)}\right)\right],\quad
x{\to}0.\label{4.14}
\end{equation}\end{La}
\begin{proof} It is immediate from \eqref{4.13} that  $E_{\beta}(x)$ is continuous in the domain $(0,\omega]$ and that
\begin{equation}  E_{\beta}(x)=\int_{0}^{1}\frac{1}{\Gamma(s)}e^{-Cs}s^{\beta-1}x^{s-1}ds+O(1).\label{4.15}
\end{equation}
Using relations
\begin{equation}\frac{1}{\Gamma(s)}=s+O(s^{2}),\quad e^{-Cs}=1+O(s),\quad s{\to}0
\label{4.16}\end{equation} we obtain
\begin{equation}  E_{\beta}(x)=\frac{1}{x}\int_{0}^{1}s^{\beta}[1+O(s)]e^{-s| \ln (x)|}ds+O(1).
\label{4.17}\end{equation}
Inequality \eqref{4.17} implies the inequality \eqref{4.14}. The lemma is proved.
\end{proof}
The next  statement \cite[p. 24]{Sakh4} follows from \eqref{4.14}. 
\begin{Pn}\label{Proposition 4.6}The operators $V_{\beta}$, defined by formula \eqref{4.11}, are bounded in  all spaces $L^{p}(0,\omega),\,p{\geq}1$.
\end{Pn}
\begin{proof} Indeed, according to \eqref{4.14} we have
\begin{equation}||V_{\beta}||_{p}{\leq}m(\beta)=
\frac{1}{\Gamma(\beta)}\int_{0}^{\omega}E_{\beta}(x)dx<\infty.\label{4.18}
\end{equation}
This proves the proposition.
\end{proof}
The operators $V_{\beta}$ have an important property:
\begin{Tm}\label{Theorem 4.7} The operators $V_{\beta}$ defined by formula \eqref{4.11}
form a semigroup, that is,
\begin{equation}V_{\alpha}V_{\beta}=V_{\alpha+\beta},\quad \alpha>0,\quad \beta>0.\label{4.19}
\end{equation}\end{Tm}
\begin{proof} The operator $V_{\alpha}V_{\beta}$ can be represented in the form
\begin{equation}V_{\alpha}V_{\beta}f=\int_{0}^{x}U_{\alpha,\beta}(x-y)f(y)dy,\label{4.20}
\end{equation} where
\begin{equation}U_{\alpha,\beta}(x)=\frac{1}{\Gamma(\alpha)\Gamma(\beta)}\int_{0}^{x}E_{\alpha}(x-y)E_{\beta}(y)dy.
\label{4.21}\end{equation}
We write the well-known formula (see \cite[Section 1.5]{BaErd}):
\begin{equation}\int_{0}^{x}(x-y)^{s-1}y^{t-1}dy=x^{s+t-1}\frac{\Gamma(s)\Gamma(t)}{\Gamma(s+t)}.
\label{4.22}\end{equation}
Using \eqref{4.13} and \eqref{4.22} we rewrite equality \eqref{4.21} in the form
\begin{equation}U_{\alpha,\beta}(x)=\frac{1}{\Gamma(\alpha)\Gamma(\beta)}\int_{0}^{\infty}
\int_{0}^{\infty}\frac{s^{\alpha-1}t^{\beta-1}}{\Gamma(s+t)}e^{-C(s+t)}x^{s+t-1}dsdt.
\label{4.23}\end{equation}
Changing the variable $s=u-t$ and the order of integration in \eqref{4.23} we have
\begin{equation}U_{\alpha,\beta}(x)=\frac{1}{\Gamma(\alpha)\Gamma(\beta)}\int_{0}^{\infty}
\frac{1}{\Gamma(u)}e^{-Cu}x^{u-1}\int_{0}^{u}(u-t)^{\alpha-1}t^{\beta-1}dtdu
\label{4.24}\end{equation}
From \eqref{4.22} and \eqref{4.24} we obtain the equality $U_{\alpha,\beta}(x)=\frac{1}{\Gamma(\alpha+\beta)}E_{\alpha+\beta}(x)$.
The theorem is proved. 
\end{proof}
2. \emph{The estimate  of $m(n\beta)$ for large values of $n$}.\\
By virtue of \eqref{4.13} and \eqref{4.18} we derive
\begin{equation}m(n\beta)=
\frac{1}{\Gamma(n\beta)}\int_{0}^{\infty}\frac{1}{s\Gamma(s)}e^{-Cs}s^{n\beta}{\omega}^{s}ds.
\label{4.25}\end{equation}
The asymptotic relations \eqref{4.5} and \eqref{4.16} imply that
\begin{equation}m(n\beta)=
\frac{1}{\Gamma(n\beta)}\left(\int_{1}^{\infty}e^{As}s^{n\beta}s^{-s}ds+O(1)\right),\label{4.26}
\end{equation} where $A=-C+1+ \ln (\omega).$ The function $e^{As}s^{n\beta}s^{-s}$
attains  maximal value at the point $s_{0}$ such that
\begin{equation}n\beta=s_0( \ln (s_0)+1-A).\label{4.27}\end{equation}
It is easy to see that
\begin{equation}\frac{n\beta}{ \ln (n\beta)}<s_0<\frac{n\gamma\beta}{ \ln (n\beta)},
\quad \gamma>1,\quad n>C_{\gamma}.\label{4.28}\end{equation}
Formula \eqref{4.28} implies that
\begin{equation}e^{As}s^{n\beta}s^{-s}<e^{|A|n\beta\gamma}
\left(\frac{n\beta\gamma}{ \ln (n\beta)}\right)^{n\beta}\left(\frac{n\beta}{ \ln (n\beta)}\right)^{-n\beta/ \ln (n\beta)}.
\label{4.29}\end{equation}Hence we have
\begin{equation}\int_{1}^{n\beta\gamma}\frac{1}{s\Gamma(s)}e^{-Cs}s^{n\beta-1}{\omega}^{s}ds
<\left(M\frac{n\beta\gamma}{ \ln (n\beta)}\right)^{n\beta}.\label{4.30}\end{equation}
The following inequalities hold for the  integrals below:
\begin{equation}\int_{n\beta\gamma}^{\infty}e^{As}s^{n\beta}s^{-s}ds<e^{An\beta}
\int_{n\beta(\gamma-1)}^{\infty}e^{-u\left( \ln (n\beta)-A\right)}du<Me^{An\beta}.\label{4.31}\end{equation}
According to \eqref{4.26}, \eqref{4.30} and \eqref{4.31}, the  estimate
\begin{equation}m(n\beta){\leq}\left(\frac{M}{ \ln (n\beta)}\right)^{n\beta} \label{4.32}\end{equation}
is valid.Using \eqref{4.32} we obtain
\begin{equation}|c_{n}(\lambda)|{\leq}m(n\beta){\leq}\left(\frac{M}{ \ln ^{\beta}(n)}\right)^{n} \label{4.33} \end{equation}
From Corollary 3.10 and relation \eqref{4.33} we have the assertion.
\begin{Pn}\label{Proposition 4.8} Let the operators  $V_{\beta}$ and $S$  be  defined by the relations \eqref{4.11} and \eqref{4.12}, respectively.
Then $\mathcal{N}(y)$, which corresponds to the operator $T=S+V_{\beta}$, satisfies the inequality 
\begin{equation} \ln \mathcal{N}(y){\leq}C\left(\frac{2M}{|y|}\right)^{1/\beta},
\quad 0<|y|{\leq}\delta\label{4.34}
\end{equation}
for some $\delta >0$.
\end{Pn}
\begin{Cy}\label{Corollary 4.9}Let the operator $T$ have the form
\begin{equation}Tf=\phi(x)f(x)+\int_{0}^{x}v(x,t)f(t)dt,\quad f{\in}H=L^{2}(0,\omega),\label{4.35}\end{equation}
where the real valued function $\phi(x)$
is bounded and  the kernel $v(x,t)$ is such that
\begin{equation}|v(x,t)|{\leq}ME_{\beta}(x-t),\quad \beta>0, \quad 0{\leq}t{\leq}x{\leq}\omega.\label{4.36}\end{equation}
Then  $\mathcal{N}(y)$, which corresponds to $T$ of the form \eqref{4.35}, satisfies \eqref{4.34}.\end{Cy}
\begin{Cy}\label{Corollary 4.10}Let the conditions of Corollary \ref{Corollary 4.9} be fulfilled
and assume that $\beta>1$. Then the operator $T$ of the form \eqref{4.35} is strongly decomposable.\end{Cy}
\begin{Rk}\label{Remark 4.11} In view of  \eqref{4.14}, estimate \eqref{4.36} in Corollaries 4.9 and 4.10 may be substituted
by simpler ones, for instance$:$
\begin{equation}
|v(x,t)| {\leq}\frac{M}{(x-t)\big(|\ln(x-t)|^{\beta+1}+1\big)}, \quad \beta>0, \quad 0{\leq}t{\leq}x{\leq}\omega. \label{4.36'}
\end{equation}
Here we used the inequality
\begin{equation}|E_{\beta}(x)|{\geq}m/(x|\ln{x}|^{\beta+1}),\quad m>0,\quad 0{\leq}x{\leq}\delta<1.
\label{4.37}\end{equation}\end{Rk}
\begin{Rk}\label{Remark 4.12}
The growth  of the resolvent
and spectral properties of the corresponding operator are closely connected.
Many interesting and fundamental results were obtained using this connection (see \cite{A, AA, AR, Lev, LyuMac, MMN, Rad}).
In   Section 4, we present important classes of such operators that their resolvents satisfy 
growth conditions essential for the study of spectral properties.
\end{Rk}
\section{Triangular integro-differential operators with difference kernels}
1. Triangular integro-differential operators with difference kernels have the form
\begin{equation}Tf=\frac{d}{dx}\int_{0}^{x}s(x-t)f(t)dt,\quad f{\in}H=L^{2}(0,\omega).\label{5.1}\end{equation}
We assume that
\begin{equation}s(x){\in}L^{2}(0,\omega),\quad x^{1/2}s^{\prime}(x){\in}L(0,\omega).\label{5.2}\end{equation}
The operators $T$ of the form \eqref{5.1} were studied in our paper \cite{Sakh3}.

Let us define the maximal chain of the orthogonal projections  \\ $E_t$ $(0\leq t \leq \omega)$:
\begin{equation}\label{op} \big(E_{t}f\big)(x)=0 \quad (0<x<\omega-t),\quad
\big(E_{t}f\big)(x)=f(x),\quad (\omega -t<x<\omega).\end{equation}
In this  section, we  prove for some classes of operators  of the form \eqref{5.1} that they  are not (S+V)-triangular
with respect to the maximal chain given by \eqref{op}.
The following example  plays an essential role in the present section.
\begin{Ee}\label{Example 5.1} Conditions \eqref{5.2} are fulfilled for the fractional integral of purely
imaginary order$:$
\begin{equation}J^{i\alpha}f=\frac{1}{\Gamma(i\alpha+1)}\frac{d}{dx}\int_{0}^{x}(x-t)^{i\alpha}f(t)dt,
\label{5.3}\end{equation} where $\alpha=\overline{\alpha}$ and $\Gamma(z)$ is the Euler function.\end{Ee}
\begin{Rk}\label{Remark 5.2}The problem of describing the spectrum of the operator $J^{i\alpha}$
was formulated in the book \cite[Ch.XXIII, Section 6]{HiPhi}.  In \cite{Sakh3}, we proved that the spectrum
of $J^{i\alpha}$ coincides with the set
\begin{equation}e^{-|\alpha|\pi/2}{\leq}|z|{\leq}e^{|\alpha|\pi/2}.\label{5.4}\end{equation} \end{Rk}
2.   Following  \cite{Sakh3} introduce the functions:
\begin{equation}\tilde{s}_{1}(\xi)=\int_{0}^{\omega}e^{it\xi}s(t)(1-t/{\omega})dt,
\quad \tilde{s}(\xi)=\int_{0}^{\omega}e^{it\xi}s(t)dt \quad\quad (\xi = \overline{\xi}). \label{5.5}\end{equation}
\begin{Dn}\label{Definition 5.3}We say that the number $\beta$ belongs to the set $\Delta^{+}\,$ if
 there is a sequence $\{\xi_k\}$ $(\xi_k=\overline{\xi_k})$ such that $\xi_{k}{\to}\infty$ and  $\big(-i\xi_{k}\tilde{s}_{1}(\xi_{k})\big){\to}\beta$ for $k \to \infty$.
 
We say that the number $\beta$ belongs to the set $\Delta^{-}\,$ if
 there is a sequence $\{\xi_k\}$ such that $\xi_{k}{\to}-\infty$ and  $\big(-i\xi_{k}\tilde{s}_{1}(\xi_{k})\big){\to}\beta$ for $k \to \infty$. 
 \end{Dn}
We shall need the next statements which were proved in  \cite{Sakh3}.
\begin{Pn}\label{Proposition 5.4}
Let a bounded operator $T$ have  the form \eqref{5.1} and let conditions
\eqref{5.2} be fulfilled.  Then
\begin{equation}||Te^{-ix\xi}-(-i\xi)e^{-ix\xi}\tilde{s}_{1}(\xi)||{\to}0,\quad
\xi{\to}\infty.\label{5.6}\end{equation}\end{Pn}
\begin{Cy}\label{Corollary 5.5} Let the conditions of Proposition \ref{Proposition 5.4} be fulfilled. Then the set $\Delta=\Delta^{+}{\bigcup}\Delta^{-}$ belongs to the spectrum of T.\end{Cy}
\begin{Tm}\label{Theorem 5.6} The operator $T$ of the form \eqref{5.1} is bounded in the space $L^{2}(0,\omega)$ if the function
$\xi{\tilde{s}(\xi)}$ is bounded on the axis $-\infty<\xi<\infty.$\end{Tm}
\begin{Pn}\label{Proposition 5.7}In the case of the operator $J^{i\alpha}$ the set $\Delta$ coincides with
the boundary of the ring \eqref{5.4}. Moreover, we have$:$

\hspace{1em}  $ z{\in}\Delta^{+}$, if $|z|=e^{-\alpha{\pi/2}}\quad$ and
$\quad z{\in}\Delta^{-}$, if $|z|=e^{\alpha{\pi/2}}.$ 
\end{Pn}
Now, we formulate our main result in this section.
\begin{Tm}\label{Theorem 5.8} Let a bounded operator $T$ have  the form \eqref{5.1} and let conditions
\eqref{5.2} be fulfilled. If the corresponding set $\Delta$ contains two different points, then the operator $T$ is not $(S+V)$-triangular
$($with respect to the maximal chain given by \eqref{op}$)$.
\end{Tm}
\begin{proof} Assuming that $T$ is $(S+V)$-triangular, that is,
$$\big(Tf\big)(x)=\phi(x)f(x)+\big(Vf\big)(x),$$ 
we prove the theorem by contradiction.
We note that  the operator $V$ is compact and $e^{ix\xi}{\to}0$ in a week sense when $|\xi|\to \infty.$ Hence, we have
$$||Te^{ix\xi}-\phi(x)e^{ix\xi}||{\to}0,\quad |\xi|{\to} \infty.$$
By comparing the last relation and \eqref{5.6} we obtain
\begin{equation} \big|(-i\xi)\tilde{s}_{1}(\xi)-\phi(x)\big|{\to}0,\quad |\xi|{\to} \infty,\label{5.7}\end{equation}
which is possible only in the case of $\Delta$ coinciding with one point.
\end{proof}
From Proposition \ref{Proposition 5.7} and Theorem \ref{Theorem 5.8}, the next
assertion easily follows.
\begin{Cy}\label{Corollary 5.9}The operator $J^{i\alpha},\, \alpha=\overline{\alpha}{\ne}0$ is not $(S+V)$-triangular
$($with respect to the maximal chain given by \eqref{op}$)$.\end{Cy}

{\bf Acknowledgements.} The author is grateful to  A. Sakhnovich and I.~Roitberg for fruitful discussions
and help in the preparation of the manuscript.

\end{document}